\documentclass{ifacconf}

\usepackage{natbib}       
\usepackage{cite}
\usepackage{amsmath,amssymb,amsfonts}
\usepackage{algorithmic}
\usepackage{graphicx}
\usepackage{algorithm,algorithmic}
\usepackage{textcomp}
\usepackage{mathtools,mathbbol,dsfont}
\usepackage{bm}
\usepackage{textcomp}
\usepackage[all]{xy}
\usepackage{chemarrow}
\usepackage{color}
\usepackage{mathrsfs,array}
\usepackage{booktabs}
\usepackage[version=4]{mhchem}

\begin{document}
\begin{frontmatter}

\title{Structure and input-to-state stability for composable computations in chemical reaction networks\thanksref{footnoteinfo}} 

\thanks[footnoteinfo]{This work was funded by the National Nature Science Foundation of
 China under Grant No. 12320101001, and the National Foreign Expert Project of China under Grant No. S20250211.}

\author[First]{Renlei Jiang,} 
\author[First,First2]{Chuanhou Gao,}
\author[Second]{Denis Dochain}

\address[First]{School of Mathematical Sciences, Zhejiang
University, Hangzhou 310058, China (e-mail: \{jiangrl,gaochou\}@zju.edu.cn).}
\address[First2]{Center for Interdisciplinary Applied Mathematics, Zhejiang University, Hangzhou, 310058, China.}
\address[Second]{ICTEAM, UCLouvain, B\^{a}timent Euler, avenue Georges Lema\^{i}tre 4-6, 1348 Louvain-la-Neuve, Belgium (e-mail: denis.dochain@uclouvain.be).}

\begin{abstract}                
In the field of molecular computation based on chemical reaction networks (CRNs), leveraging parallelism to enable coupled mass-action systems (MASs) to retain predefined computational functionality has been a research focus. MASs exhibiting this property are termed composable. This paper investigates the structural conditions under which two MASs are composable. By leveraging input-to-state stability (ISS) property, we identify a specific class of CRN architectures that guarantee composability with other networks. A concrete example demonstrates the validity of this conclusion and illustrates the application of composability in computing composite functions.
\end{abstract}


\end{frontmatter}

\section{Introduction}
In recent years, the programmability of biochemical systems through molecular computation has attracted substantial research interest, representing a convergence of multiple disciplines including computer science, mathematics, and systems biology. This emerging field uses molecular components, like DNA and enzymes, as computational elements to design systems that can process information, make decisions, and execute complex behaviors in biological environments. Unlike conventional silicon-based hardware computation, molecular computation operates under energy-efficient, parallelized, and biocompatible conditions, making it particularly suitable for applications such as medical diagnosis \citep{zhang2020cancer} and data storage \citep{wang2023parallel}.

Chemical reaction networks (CRNs) serve as a primary framework for molecular computation. They characterize interactions between chemical species through reactions and are widely used to represent complex biochemical processes such as metabolism and gene regulation. Mathematically, CRNs governed by mass-action kinetics, referred to as mass-action systems (MASs), can be modeled by a group of ordinary differential equations (ODEs), which describe the time evolution of species concentrations. The structure and parameters (reaction rate constants) of a MAS determine its dynamical behavior, including properties such as stability \citep{feinberg1987chemical}, persistence \citep{craciun2013persistence} and oscillations \citep{domijan2009bistability}. This versatility renders CRNs powerful tools for understanding both natural and engineered chemical systems, while also providing a robust foundation for molecular computation.

As studies have demonstrated that any MAS can be physically implemented via DNA strand displacement reactions \citep{soloveichik2010dna}, a growing number of researchers began to dedicate efforts to designing biochemical systems capable of complex information processing based on MAS frameworks \citep{anderson2021reaction,chen2024synthetic,fan2025automatic_1}. However, these studies continue to face a fundamental challenge: an inherent conflict exists between the parallel nature of chemical reactions and the sequential execution required of computation. One feasible approach to address this issue is to use chemical oscillators \citep{shi2025controlling}, which exhibit phases where only a single species concentration is strictly greater than zero. This characteristic enables the selective catalysis of specific reactions while simultaneously suppressing others, thereby segregating two or more inherently parallel reactions and orchestrating their sequential execution. But this approach tends to render the network overly complex and introduces significant errors when applied to complex multi-step molecular computation tasks \citep{fan2025automatic_2}. Consequently, it is imperative to investigate strategies that harness the inherent parallelism of chemical reactions for molecular computations, enabling the concurrent execution of multiple reactions to achieve desired computational results. \citep{chalk2019composable} formalized this concept as composability, which means the computation of a composite function can be achieved by interconnecting CRNs that compute elementary functions. They investigated the composability of deterministic rate-independent CRNs and precisely characterized the class of functions computable by such CRNs. \citep{severson2019composable} extended this line of research to stochastic rate-independent CRNs. \citep{jiang2025input} investigated the composability of MAS and developed a novel method to determine composability based on input-to-state stability (ISS) criteria. In contrast to the severely constrained computational capacity of rate-independent CRNs \citep{chen2023rate}, MASs have been proven to be Turing-universal \citep{fages2017strong}, implying that any computation can be embedded into a class of polynomial ODEs and implemented by CRNs. Consequently, research on the composability of MAS facilitates the realization of a broader spectrum of computational functions.

Based on the definitions and conclusions established in \citep{jiang2025input}, this paper establishes a critical connection between network structure and composability, investigating the structural conditions required for MAS to achieve composable computation. Our paper is organized as follows. Some necessary preliminaries are given in Section \ref{section_2}. Section \ref{section_3} is dedicated to establish the network structural conditions for composable MAS-based molecular computations. In Section \ref{section_4}, we present an example to illustrate our theoretical conclusions, followed by some discussions about our work.

\noindent{\textbf{Mathematical Notation:}}\\
\rule[1ex]{\columnwidth}{0.8pt}
\begin{description}
\item[\hspace{0em}{$\mathbb{R}^n, \mathbb{R}^n_{\geq 0},\mathbb{R}^n_{>0},\mathbb{Z}^n_{\ge 0}$}]: $n$-dimensional real space; $n$-dimensional non-negative real space; $n$-dimensional positive real space; $n$-dimensional non-negative integer space.
\item[\hspace{0em}{$s^{v_{\cdot j}}$}]: $s^{v_{\cdot j}}\triangleq\prod_{i=1}^{n}s_{i}^{v_{ij}}$, where $s,v_{\cdot j}\in\mathbb{R}^{n}$.
\item[\hspace{0em}{$\mathcal{K},\mathcal{K}_{\infty},\mathcal{KL}$}]: $\gamma : \mathbb{R}_{\ge 0} \to \mathbb{R}_{\ge 0}$ is a class $\mathcal{K}$ function if it is continuous strictly increasing and satisfies $\gamma(0)=0$; $\gamma$ is a class $\mathcal{K}_\infty$ function if it is a class $\mathcal{K}$ function and satisfies $\lim_{s \to \infty}\gamma(s)=\infty$; $\beta : \mathbb{R}_{\ge 0} \times \mathbb{R}_{\ge 0} \to \mathbb{R}_{\ge 0}$ is a class $\mathcal{KL}$ function if for each fixed $t$ the mapping $\beta(\cdot,t)$ is a class $\mathcal{K}$ function and for each fixed $s$ it decreases to zero on $t$ as $t \to \infty$.
\end{description}
\rule[1ex]{\columnwidth}{0.8pt}


\section{Preliminaries}\label{section_2}

In this section, we introduce some basic concepts about CRNs.

Consider a CRN with $n$ species, denoted by $S_1,...,S_n$, and $r$ reactions with the $j$th reaction written as
$$\sum_{i=1}^nv_{ij}S_i\to \sum_{i=1}^nv'_{ij}S_i,$$
where $v_{.j}, v'_{.j}\in\mathbb{Z}_{\geq 0}^n$ represent the reactant complex and the product complex of the reaction, respectively. For simplicity, this reaction is often written as $v_{.j}\to v'_{.j}$. We often use triple $(\mathcal{S},\mathcal{C},\mathcal{R})$ to express a CRN, where $\mathcal{S}$ is \textit{species set}, $\mathcal{C}$ is \textit{complex set}, and $\mathcal{R}$ is \textit{reaction set}. The \textit{stoichiometric subspace} of a CRN is defined by $\mathscr{S}=\textrm{span}\{ v_{\cdot 1}^{\prime}-v_{\cdot 1},...,v_{\cdot r}^{\prime}-v_{\cdot r}\}.$ For $s_0 \in \mathbb{R}^{n}_{\ge 0}$, the set $s_0+\mathscr{S}=\{s_0+s|s \in \mathscr{S}\}$, $\left(s_0+\mathscr{S}\right)\cap \mathbb{R}^n_{\ge 0}$ and $\left(s_0+\mathscr{S}\right)\cap \mathbb{R}^n_{> 0}$ are called the \textit{stoichiometric compatibility class}, the \textit{nonnegative stoichiometric compatibility class}, and the \textit{positive stoichiometric compatibility class} of $s_0$, respectively.

When a CRN is equipped with mass-action kinetics, the rate of reaction $v_{\cdot j} \to v_{\cdot j}^{\prime}$ is measured by $k_js^{v_{\cdot j}}$, where $k_j>0$ represents the rate constant, and $s \in \mathbb{R}_{\ge 0}^n$ with each element $s_i~(i=1,...,n)$ to represent the concentration of the species $S_i$. The dynamics of a \textit{MAS} that describes the change of concentrations of all species over time $t$ is given by
\begin{equation}\label{general dynamics}
    \frac{\mathrm{d}s(t)}{\mathrm{d}t}=\sum_{j=1}^{r}k_js^{v_{\cdot j}}\left(v_{\cdot j}^{\prime}-v_{\cdot j} \right).
\end{equation}
We often use a quadruple $(\mathcal{S},\mathcal{C},\mathcal{R},\kappa)$ to express a MAS. If $\kappa$ is not always constant but time-varying, we call $\left(\mathcal{S},\mathcal{C},\mathcal{R},\kappa(t)\right)$ a \textit{generalized MAS}. A constant vector $\bar{s} \in \mathbb{R}^n_{>0}$ is called a \textit{positive equilibrium} of system (\ref{general dynamics}) if 
\begin{equation}
    \sum_{j=1}^{r}k_j\bar{s}^{v_{\cdot j}}\left(v_{\cdot j}^{\prime}-v_{\cdot j} \right)=0.
\end{equation}
By integrating (\ref{general dynamics}) from $0$ to $t$, we get
\begin{equation}\label{integ_dynamics}
 s(t) = s_0 + \sum_{j=1}^r \left(v_{\cdot j}^{\prime}-v_{\cdot j} \right) \int_0^t k_js^{v_{\cdot j}}(\tau)\mathrm{d}\tau,   
\end{equation}
where $s_0=s(0)$ is the initial state of $(\mathcal{S},\mathcal{C},\mathcal{R},\kappa)$. It is clear that the state of $(\mathcal{S},\mathcal{C},\mathcal{R},\kappa)$ will evolve in the nonnegative stoichiometric compatibility class of $s_0$, i.e., in $(s_0+\mathscr{S}) \bigcap \mathbb{R}_{\ge 0}^n$, where $\mathscr{S}=\text{span}\{v_{\cdot 1}^{\prime}-v_{\cdot 1},...,v_{\cdot r}^{\prime}-v_{\cdot r}\}$ We say the trajectory $s(t)$ is \textit{persistent} if 
$$\liminf_{t \to \infty}s_i(t)>0,\quad i=1,2,\cdots ,n.$$

The structure of CRNs has been a focal point in research, as CRNs with specific architectures can exhibit particular dynamical behaviors in their corresponding MAS (regardless of rate constants values $\kappa$). A CRN can be seen as a directed graph, where the nodes represent complexes and the directed edges correspond to reactions. Each connected component of the directed graph is called a \textit{linkage class}. A CRN is \textit{weakly reversible} if each of its linkage classes is strongly connected. A CRN is \textit{reversible} if the reactions in $\mathcal{R}$ are symmetric, that is, if $\forall v_{\cdot j} \to v_{\cdot j}^{\prime} \in \mathcal{R}$ it holds $v_{\cdot j}^{\prime} \to v_{\cdot j} \in \mathcal{R}$. A reversible CRN is necessarily weakly reversible, whereas the converse does not hold. For a CRN, let $n$ denote the number of complexes and $l$ denote the number of linkage classes. The \textit{deficiency} of the CRN is defined by $\delta =n-l-\dim \mathscr{S}$. It is noted that the deficiency of a CRN is nonnegative because it can be seen as the dimension of a certain linear subspace \citep{feinberg2019foundations}.

\begin{exmp}\label{ex_deficiency}
    Consider the following CRN
    \begin{equation*}
        Z_1 \ce{<=>[][]} Z_2,
    \end{equation*}
    which is reversible (surely weakly reversible). In addition, it has two complexes $(1,0)^\top$ and $(0,1)^\top$, one linkage class, and its stoichiometric subspace is $\mathscr{S}=\text{span}\{(1,-1)^{\top},(-1,1)^{\top}\}$. Hence, its deficiency is $\delta=2-1-1=0$. When the CRN is endowed with mass-action kinetics $\kappa=\{1,2\}$, forming a MAS, its dynamics are governed by
    \begin{equation}\label{ex_MAS}
    \begin{cases}
        \dot{z}_1=2z_2-z_1, \\
        \dot{z}_2=z_1-2z_2.
    \end{cases}
    \end{equation}
    Similarly, when equipped with generalized mass-action kinetics $\kappa=\left(k_1(t),k_2(t)\right)$, it constitutes a generalized MAS whose dynamics are described by
    \begin{equation}
    \begin{cases}
        \dot{z}_1=k_2(t)z_2-k_1(t)z_1, \\
        \dot{z}_2=k_1(t)z_1-k_2(t)z_2.
    \end{cases}
    \end{equation}
\end{exmp}

One of the most famous results in the study of CRN structure is the \textit{Deficiency Zero Theorem} \citep{feinberg1987chemical}, which completely characterizes the dynamical behavior of MAS with zero deficiency. It tell us that if a MAS $(\mathcal{S},\mathcal{C},\mathcal{R},\kappa)$ is weakly reversible and has zero deficiency, then for any $\kappa \in \mathbb{R}^{r}_{>0}$, the dynamics of MAS has the following properties: There exists within each positive stoichiometric compatibility class precisely one equilibrium $\bar{s}=\bar{s}(s_0,\kappa)$, and the pseudo-Helmholtz free energy function
\begin{equation}\label{Lyafun_deficiency_zero}
    V(s,\bar{s})=\sum_{j=1}^r \left(s_j(\ln s_j-\ln \bar{s}_j-1)-\bar{s}_j\right),~s \in \mathbb{R}^{n}_{> 0}
\end{equation}
is a Lyapunov function that ensures $\bar{s}$ to be locally asymptotically stable. Since the CRN in \textit{Example} \ref{ex_deficiency} is weakly reversible and has zero deficiency, the dynamics (\ref{ex_MAS}) has the unique locally asymptotically stable equilibrium $(\bar{z}_1,\bar{z}_2)=\left(2C/3,C/3 \right)$ in the positive stoichiometric compatibility class $\{(z_1,z_2) \in \mathbb{R}^2_{>0}|z_1+z_2=C\}$ with $C>0$.


\section{Structure condition for composable msCRCs}\label{section_3}

\subsection{MAS-based composable molecular computations}
In the study of molecular computation, MAS serves as a prevalent mathematical framework designed to represent computational inputs via species concentrations and to encode computational outputs as the limiting steady state of the relevant species concentrations. Jiang et. al. (2025) established a framework for MAS-based molecular computation and investigated strategy for achieving their composite computations.

\begin{defn}[msCRC]
    A mass-action chemical reaction computer (msCRC) is a tuple $\mathscr{C}=(\mathcal{S},\mathcal{C},\mathcal{R},\kappa,\mathcal{X},\mathcal{Y})$, where $(\mathcal{S,C,R},\kappa)$ is a MAS, $\mathcal{X} \subset \mathcal{S}$ is the input species set, and $\mathcal{Y}=\mathcal{S} \setminus \mathcal{X}$ is the output species set.
\end{defn}

\begin{defn}[Dynamic computation]\label{Def_dynamic_compute} Given a msCRC $\mathscr{C}=(\mathcal{S},\mathcal{C},\mathcal{R},\kappa,\mathcal{X},\mathcal{Y})$ with $m~(m<n)$ species in $\mathcal{X}$, and a positive function $\sigma:\mathbb{R}^{m}_{\ge 0} \to \mathbb{R}^{n-m}_{\ge 0}$, denote the dynamics of this msCRC by
    \begin{equation}\label{CRC_dynamics}
    \begin{cases}
        \dot{x}=f(x,y), \\
        \dot{y}=g(x,y),
    \end{cases} x(0)=x_0,y(0)=y_0,
    \end{equation}
where $x \in \mathbb{R}^{m}_{\ge 0}$ and $y \in \mathbb{R}^{n-m}_{\ge 0}$ are the concentrations of input species and output species, respectively. The msCRC is a dynamic computation of function $\sigma$, if the output species concentrations satisfy
    \begin{equation}\label{eq:dyCom}
  \lim_{t \to \infty}y(t)=\sigma(x_0).      
    \end{equation}
\end{defn}

For MAS-based molecular computation, a critical issue involves investigating the composability of multiple simple molecular computing systems. Specifically, it is essential to explore the possibility of leveraging their combinations to achieve the corresponding composite computations.


\begin{defn}[Dynamically Composable]\label{Def_dynamic_composable}
Given two msCRCs $\mathscr{C}^\mathtt{1}=(\mathcal{S}^\mathtt{1},\mathcal{C}^\mathtt{1},\mathcal{R}^\mathtt{1},\kappa^\mathtt{1},\mathcal{X}^\mathtt{1},\mathcal{Y}^\mathtt{1})$ and $\mathscr{C}^\mathtt{2}=(\mathcal{S}^{\mathtt{2}},\mathcal{C}^{\mathtt{2}},\mathcal{R}^{\mathtt{2}},\kappa^{\mathtt{2}},\mathcal{X}^{\mathtt{2}},\mathcal{Y}^{\mathtt{2}})$ that satisfy the following assumptions:
\begin{enumerate}
    \item[(A.1)] $\mathcal{Y}^\mathtt{1}=\mathcal{X}^{\mathtt{2}},~\mathcal{Y}^{\mathtt{2}}\cap \mathcal{X}^{\mathtt{1}}=\varnothing$;
    \item[(A.2)] $\mathscr{C}^\mathtt{1}$ and $\mathscr{C}^\mathtt{2}$ have respective dynamics to be
    \begin{equation}\label{dx=f}
        \begin{cases}
          \dot{x}^\mathtt{1}=f^\mathtt{1}(x^\mathtt{1},y^\mathtt{1}), \\
          \dot{y}^\mathtt{1}=g^\mathtt{1}(x^\mathtt{1},y^\mathtt{1}),
        \end{cases}x^\mathtt{1}(0)=x^\mathtt{1}_0,~y^\mathtt{1}(0)=y^\mathtt{1}_0
    \end{equation} 
    and
    \begin{equation}\label{dy=g}
        \begin{cases}
            \dot{x}^\mathtt{2}=0, \\
            \dot{y}^\mathtt{2}=g^\mathtt{2}(x^\mathtt{2},y^\mathtt{2}),
        \end{cases} x^\mathtt{2}(0)=\bar{y}^\mathtt{1},~y^\mathtt{2}(0)=y^\mathtt{2}_0;
    \end{equation}
    \item[(A.3)] $\lim_{t \to \infty} y^\mathtt{1}(t)= \bar{y}^\mathtt{1},~\lim_{t \to \infty} y^{\mathtt{2}}(t)= \bar{y}^{\mathtt{2}}$.
\end{enumerate}
If the solution of the coupled system 
\begin{equation}\label{coupled_system}
    \begin{cases}
     \dot{x}^\mathtt{1}=f^\mathtt{1}(x^\mathtt{1},y^\mathtt{1}), \\
     \dot{y}^\mathtt{1}=g^\mathtt{1}(x^\mathtt{1},y^\mathtt{1}),\\
     \dot{y}^\mathtt{2}=g^\mathtt{2}(x^\mathtt{2},y^\mathtt{2}), 
    \end{cases} \begin{array}{l}
         x^\mathtt{1}(0)=x^\mathtt{1}_0,~y^\mathtt{1}(0)=y^\mathtt{1}_0, \\
         y^\mathtt{2}(0)=y^\mathtt{2}_0,
    \end{array}
\end{equation}
satisfies 
\begin{equation}\label{dynamic_composable}
    \lim_{t \rightarrow \infty} \left [ \begin{array}{c}
         y^\mathtt{1}(t)  \\
         y^\mathtt{2}(t)  \\
    \end{array} \right] =\left[ \begin{array}{c}
         \bar{y}^\mathtt{1} \\
         \bar{y}^\mathtt{2} \\
    \end{array} \right],
\end{equation}
then the two msCRCs are said to be dynamically composable. 
\end{defn}

Assumptions (A.1)-(A.3) are crucial for studying the composability of msCRCs, which will be repeatedly referenced in the following theoretical developments. If $\mathscr{C}^{\mathtt{1}}$ and $\mathscr{C}^{\mathtt{2}}$ are regarded as two msCRCs implementing distinct molecular computations, assumptions (A.1)-(A.2) imply that the output of $\mathscr{C}^{\mathtt{1}}$ serves as the input of $\mathscr{C}^{\mathtt{2}}$, while the species in $\mathscr{C}^{\mathtt{1}}$ remain unchanged in $\mathscr{C}^{\mathtt{2}}$ (i.e., they act as catalysts). Assumption (A.3) further indicates that $\mathscr{C}^{\mathtt{1}}$ and $\mathscr{C}^{\mathtt{2}}$ can achieve layer-by-layer computation of a composite function. The limiting steady state of the coupled system (\ref{coupled_system}) can be thought as a layer-by-layer composition of their individual output limiting steady state.

\begin{thm}\citep{jiang2025input}\label{thm_ISS}
    Suppose $\mathscr{C}^\mathtt{1}$ and $\mathscr{C}^\mathtt{2}$ satisfy assumptions (A.1)-(A.3). Consider the $y^\mathtt{2}$-related part of the system (\ref{dy=g}), if it is ISS regarding $(\bar{y}^{\mathtt{1}},\bar{y}^{\mathtt{2}})$, that is, there exist a class $\mathcal{KL}$ function $\beta$ and a class $\mathcal{K}$ function $\gamma$ such that for any $y^{\mathtt{2}}_0$ and any bounded input $x^{\mathtt{2}}(t)$, the solution $y^{\mathtt{2}}(t)$ exists for all $t \ge 0$ and satisfies
        \begin{equation}\label{ISS_def}
           | y^{\mathtt{2}} \left( t \right) -\bar{y}^{\mathtt{2}}| \le \beta \left( | y^{\mathtt{2}}_0 - \bar{y}^{\mathtt{2}} | , t \right) + \gamma ( \sup \limits_{0\le \tau \le t}| x^{\mathtt{2}} \left( \tau \right) -\bar{y}^{\mathtt{1}} | ),
        \end{equation}   
    then $\mathscr{C}^\mathtt{1}$ and $\mathscr{C}^\mathtt{2}$ are dynamically composable.
\end{thm}

\textit{Theorem} \ref{thm_ISS} shows that (\ref{ISS_def}) is a sufficient condition for composability. However, it is not easy to find functions $\beta$ and $\gamma$ satisfying (\ref{ISS_def}). In the following subsection, we will investigate sufficient conditions based on the structure of MAS that admit more straightforward verifications.

\subsection{Structure condition for composable msCRCs}

The application of \textit{Theorem} \ref{thm_ISS} requires ISS conditions. A standard method for proving ISS is to construct an appropriate Lyapunov function.

\begin{lem}\label{lem1}
    Suppose $\lim_{t \to \infty} u(t) = \bar{u}$, and $s(t)$ is a bounded solution to the generalized MAS
    \begin{equation}\label{generalized_MAS}
        \dot{s}=f(u,s)=\sum_{j=1}^{r}u_j(t)s^{v_{\cdot j}}\left(v_{\cdot j}^{\prime}-v_{\cdot j}\right),~ s(0)=s_0,
    \end{equation}
    where $s(t) \in \mathbb{R}^n_{\ge 0},~u(t):[0,+\infty) \to \mathbb{U}$ is the input function with $\mathbb{U} \subset \mathbb{R}^r_{> 0}$. Let $\bar{s} \in (s_0+\mathscr{S})\cap \mathbb{R}^n_{> 0}$ be the unique equilibrium of the MAS $\dot{s}=f(\bar{u},s)$ in the positive stoichiometric compatibility class. In addition, suppose that there is a function $V: \mathbb{R}^n_{\ge 0} \to \mathbb{R}_{\ge 0}$ that satisfies the following:
    \begin{enumerate}
        \item the restriction of $V$ to $\mathbb{R}^n_{>0}$ is continuously differentiable;
        \item $V(s) \to +\infty$ as $|s| \to +\infty$;
        \item $V(\bar{s})=0$ and $V(s) > 0,~s \in \mathbb{R}^n_{\ge 0} \backslash \{\bar{s}\}$;
        \item for each compact set $F \subset \mathbb{R}^n_{\ge 0}$, there exist class $\mathcal{K}_{\infty}$ functions $\alpha,\rho$ such that $$dV/dt=\nabla V^{\top}(s)f(u,s) \le -\alpha(|s-\bar{s}|)+\rho (|u-\bar{u}|)$$
        for all $u \in \mathbb{U}$ and $s \in F \cap (s_0+\mathscr{S})\cap \mathbb{R}^n_{>0},$
    \end{enumerate}
    then either $s(t) \to \bar{s}$ or $s(t) \to \partial \mathbb{R}^n_{\ge 0}$.
\end{lem}

The function $V(s)$ in \textit{Lemma} \ref{lem1} is called an ISS-Lyapunov function \citep{chaves2002state}, which is only required to be differentiable in the positive orthant. In this regard, it is slightly different from the traditional ISS-Lyapunov function in \citep{sontag1995characterizations}, which implies the ISS property and consequently ensures $s \to \bar{s}$ (see \citep{jiang2025input}). The reason why \textit{Lemma} \ref{lem1} fails to guarantee $s \to \bar{s}$ lies in this slight difference. Therefore, to ensure $s \to \bar{s}$ under this condition, the exclusion of another case is required.

\begin{lem}\label{lem2}
    Under the condition of \textit{Lemma} \ref{lem1}, if $s(t)$ is persistent, then $s(t) \to \bar{s}$.
\end{lem}

Before introducing the main theorem of this section, an additional concept requires exposition.
\begin{defn}[Reduced system]
    Suppose $\mathscr{C}=(\mathcal{S},\mathcal{C},\mathcal{R},\kappa, \\ \mathcal{X},\mathcal{Y})$ is a msCRC. A generalized MAS $\tilde{\mathscr{C}}=(\tilde{\mathcal{S}},\tilde{\mathcal{C}},\tilde{\mathcal{R}},\tilde{\kappa}(t))$ is called a reduced system of $\mathscr{C}$ if
    \begin{enumerate}
        \item $\tilde{\mathcal{S}}=\mathcal{Y}$;
        \item $\tilde{\mathcal{C}}=\pi_{\mathcal{Y}}(\mathcal{C})$, where $\pi_{\mathcal{Y}}$ denotes the projection of each complex $v_{\cdot j} \in \mathcal{C}$ onto species in $\mathcal{Y}$;
        \item $\tilde{\mathcal{R}}=\{\pi_{\mathcal{Y}}(v_{\cdot j})\overset{\tilde{\kappa}_j(t)}{\longrightarrow} \pi_{\mathcal{Y}}(v_{\cdot j}^{\prime}),~j=1,...,r\}$, where the reaction rate is given by $\tilde{\kappa}_j(t)=\kappa_j \Pi_{S_i \notin \tilde{S}}s_i^{v_{ij}}(t)$.
    \end{enumerate}
\end{defn}

\begin{exmp}\label{ex_reduced}
    Consider a msCRC
    \begin{equation}
        X \overset{k_1}{\longrightarrow}2X, \quad X+Y \overset{k_2}{\longrightarrow}2Y, \quad Y \overset{k_3}{\longrightarrow}\varnothing
    \end{equation}
    with $\mathcal{X}=\{X\},~\mathcal{Y}=\{Y\}$, which has dynamics
    \begin{equation*}
    \begin{cases}
        \dot{x}=k_1x-k_2xy, \\
        \dot{y}=k_2xy-k_3y.
    \end{cases}
    \end{equation*}
    Its reduced system is
    \begin{equation}
        \varnothing \overset{k_3}{\longleftarrow} Y \overset{k_2x(t)}{\longrightarrow}2Y
    \end{equation}
    with dynamics to be $\dot{y}=k_2x(t)y-k_3y$.
\end{exmp}

The reduced system is constructed by incorporating the dynamics of the input species into the rate constants, thereby removing the input species to achieve system simplification. This concept is widely used in the theoretical analysis of CRN \citep{anderson2011proof,zhang2025network}. 
It is noted that the reduced system preserves the same dynamics as the output species of the original system. This property enables us to analyze the dynamical behavior of the original system when coupled with other systems through the study of its reduced system. Based on this we give our main theorem, which enables direct determination of composability from the structural properties of the MAS.

\begin{thm}\label{thm_main}
    Suppose that $\mathscr{C}^\mathtt{1}$ and $\mathscr{C}^\mathtt{2}$ satisfy the assumptions (A.1)-(A.3) of \textit{Definition} \ref{Def_dynamic_composable}. Let $\tilde{\mathscr{C}}^{\mathtt{2}}$ be the reduced system of $\mathscr{C}^{\mathtt{2}}$. If $\tilde{\mathscr{C}}^{\mathtt{2}}$ satisfies
    \begin{enumerate}
        \item it is weakly reversible;
        \item it has single linkage class and zero deficiency;
        \item it is mass-conservative, that is, $\exists v \in \mathbb{R}^n_{>0}$ such that
        $$v^\top g^{\mathtt{2}}(x^{\mathtt{2}},y^{\mathtt{2}})=0,~\forall x^{\mathtt{2}} \in \mathbb{R}^{n_2-m_2}_{>0}, \forall y^{\mathtt{2}} \in \mathbb{R}^{m_2}_{> 0},$$
    \end{enumerate}
    then $\mathscr{C}^\mathtt{1}$ and $\mathscr{C}^\mathtt{2}$ are dynamically composable.
\end{thm}

\begin{rem}
    The proof of \textit{Theorem} \ref{thm_main} proceeds in two stages. First, we indicate that (\ref{Lyafun_deficiency_zero}) indeed constitutes an ISS-Lyapunov function for the generalized MAS $\tilde{\mathscr{C}}^{\mathtt{2}}$, which is a previous result in \citep{chaves2005input}. Second, we employ this function to establish composability, which remains previously unaddressed.
\end{rem}

All conditions in \textit{Theorem} \ref{thm_main} can be directly observed and verified from the structures of $\mathscr{C}^{\mathtt{1}}$ and $\mathscr{C}^{\mathtt{2}}$. This facilitates the deduction and understanding of composability through the structural analysis of CRNs.


\section{Example and Discussions}\label{section_4}


To get a deeper understanding of composable msCRCs with specific architectures and their applications in molecular computations, we provide a concrete example of a molecular computing system.
\begin{exmp}\label{ex}
    Consider a msCRC
    \begin{equation}\label{CRN_1}
    \begin{split}
        & X_1 \overset{1}{\longrightarrow} X_1+Y_1, \quad X_2 \overset{1}{\longrightarrow} X_2+Y_1, \quad Y_1 \overset{1}{\longrightarrow} \varnothing, \\
        & X_3 \overset{1}{\longrightarrow} X_3+Y_2, \quad X_4 \overset{1}{\longrightarrow} X_4+Y_2, \quad Y_2 \overset{1}{\longrightarrow} \varnothing,
    \end{split}
    \end{equation}
\end{exmp}
with $\mathcal{X}_1=\{X_1,X_2,X_3,X_4\},~\mathcal{Y}_1=\{Y_1,Y_2\}$. It has dynamics and solution to be 
\begin{equation*}
\begin{split}
    &\begin{cases}
        \dot{x}_i=0, \\
        \dot{y}_j=x_{2j-1}+x_{2j}-y_j,
    \end{cases} \\
    \Rightarrow &\begin{cases}
        x_i(t)=x_i(0), \\
        y_j(t)=x_{2j-1}(0)+x_{2j}(0)+c_je^{-t},
    \end{cases}
\end{split}
\end{equation*}
where $c_j=y_j(0)-x_{2j-1}(0)-x_{2j}(0),~i=1,2,3,4,~j=1,2.$ Then we have $\lim_{t \to \infty}y_j(t)=x_{2j-1}(0)+x_{2j}(0),~j=1,2.$ This means that (\ref{CRN_1}) can dynamically compute the function $\sigma_1: \mathbb{R}^4_{\ge 0} \to \mathbb{R}^2_{\ge 0}$ given by
$$\sigma_1(x_1,x_2,x_3,x_4)=\left(x_1+x_2,x_3+x_4 \right)^\top.$$ 
Consider another msCRC 
\begin{equation}\label{CRN_2}
        Y_1+Z_1 \overset{1}{\longrightarrow}Y_1+Z_2, \quad  Y_2+Z_2 \overset{1}{\longrightarrow}Y_2+Z_1,
    \end{equation}
    with $\mathcal{X}=\{Y_1,Y_2\},~\mathcal{Y}=\{Z_1,Z_2\}$, which has dynamics
    \begin{equation}\label{eq_CRN_2}
    \begin{cases}
        \dot{y}_1=\dot{y}_2=0, \\
        \dot{z}_1=y_2z_2-y_1z_1, \\
        \dot{z}_2=y_1z_1-y_2z_2.
    \end{cases}
    \end{equation}
At the limiting steady state, we have $\bar{z}_1/\bar{z}_2=y_2/y_1=y_2(0)/y_1(0)$, and $\bar{z}_1+\bar{z}_2=z_1(0)+z_2(0)$. Therefore, by setting $z_1(0)+z_2(0)=1$ we can get $\lim_{t \to \infty}z_j(t)=\bar{z}_j=\frac{y_j(0)}{y_1(0)+y_2(0)},~j=1,2.$ This shows that (\ref{CRN_2}) is a dynamic computation of the \textit{Normalization Function} $\sigma_2: \mathbb{R}^2_{\ge 0} \to \mathbb{R}^2_{\ge 0}$ with
$$\sigma_2(y_1,y_2)=\left(\frac{y_1}{y_1+y_2},\frac{y_2}{y_1+y_2} \right)^{\top}.
$$
Note that the reduced system of (\ref{CRN_2}) is
    \begin{equation}\label{CRN_2_reduced}
        Z_1 \ce{<=>[$y_1(t)$][$y_2(t)$]} Z_2,
    \end{equation}
by \textit{Example} \ref{ex_deficiency} it is weakly reversible and has single linkage class and zero deficiency. In addition, by (\ref{eq_CRN_2}) we have $g^{\mathtt{2}}(y,z)=(y_2z_2-y_1z_1,y_1z_1-y_2z_2)^{\top}$. Since $v=(1,1) \in \mathbb{R}^{2}_{> 0}$ satisfies $v^{\top}g^{2}=0$, we know that (\ref{CRN_2_reduced}) is mass-conservative. Thus, by \textit{Theorem} \ref{thm_main} we know that (\ref{CRN_1}) and (\ref{CRN_2}) are dynamically composable.

\begin{figure}[ht]
    \centering
    \includegraphics[width=0.47\textwidth]{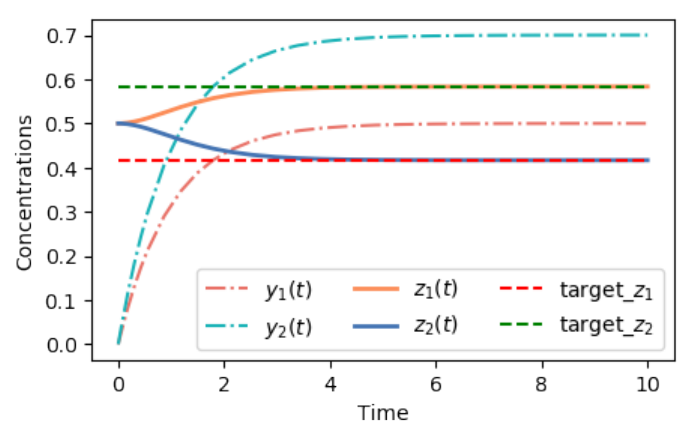}
    \caption{The computation results of the function (\ref{fun_composite})} with initial values $x_1(0)=0.2,x_2(0)=0.3,x_3(0)=0.6,x_4(0)=0.1,~y_1(0)=y_2(0)=0,~z_1(0)=z_2(0)=0.5$.
    \label{fig1}
\end{figure}

\textit{Example} \ref{ex} enables direct determination of composability from structural information of the two networks, with no need for analyzing dynamics of networks and using \textit{Theorem} \ref{thm_ISS}. This composability implies that the composition of the two networks can serve to compute the composite function $\sigma_2 \circ \sigma_1$ (see \textit{Theorem 3.1} in \citep{jiang2025input}), that is, compute the function
\begin{equation}\label{fun_composite}
    \sigma(x_1,x_2,x_3,x_4)=\left(\frac{x_1+x_2}{\sum_{i=1}^4x_i},\frac{x_3+x_4}{\sum_{i=1}^4x_i} \right).
\end{equation}
Fig. \ref{fig1} shows the simulations for the composition of (\ref{CRN_1}) and (\ref{CRN_2}) that serve to compute (\ref{fun_composite}). It should be noted that in \textit{Example} \ref{ex}, the computation of $\sigma_1$ permits arbitrary initial values for $y_1$ and $y_2$ (here, $y_1(0)=y_2(0)=0$), whereas the computation of $\sigma_2$ requires the constraint $z_1(0)+z_2(0)=1$ to be satisfied (here, $z_1(0)=z_2(0)=0.5$).

The zero deficiency network (\ref{CRN_2_reduced}) in \textit{Example} \ref{ex} has only two species and hence is relatively simple and exhibits straightforward dynamics. However, some zero deficiency networks may comprise hundreds of species and reactions, rendering the analysis of their dynamics, particularly their dynamics when composited with other networks, exceptionally complex. In such cases, the advantage of \textit{Theorem \ref{thm_main}}, which eliminates the need to analyze the dynamics of the composited network, becomes particularly pronounced. Future work will continue to focus on the influence of network structure on composability, exploring a wider range of composable network architectures and their computable functions, thereby constructing a "composable elementary msCRC library".

\bibliography{ifacconf}             

\appendix
\section{Proofs of lemmas and theorem}  
In this appendix, we provide all proofs to the results given in the previous sections.

\textbf{Proof of Lemma \ref{lem1}:} Suppose that $s(t) \nrightarrow \bar{s}$, then there exist a $\eta_0>0$ and a sequence $\{t_k,k\ge 1\}$ satisfying $t_k \to \infty$ and $|s(t_k)-\bar{s}|\ge \eta_0$. Since $u(t) \to \bar{u}$, for $\xi=\alpha(\eta_0)/2>0$, there exists a $T>0$ such that $\rho(|u(t)-\bar{u}|)<\xi$. Then for any $t>T$ and $|s(t)-\bar{s}|\ge \eta_0$, we have
    \begin{equation}\label{pf2}
        \begin{split}
            dV/dt &\le  -\alpha(|s(t)-\bar{s}|)+\rho(|u(t)-\bar{u}|) \\
            & \le -\alpha(\eta_0)+\xi \\
            & = -\alpha(\eta_0)/2,
        \end{split}
    \end{equation}
    where the first inequality holds due to the definition of ISS-Lyapunov function and $s(t)$ being bounded. If $\exists ~t_0>T$ such that $|s(t_0)-\bar{s}|<\eta_0$, by (\ref{pf2}) we know $|s(t)-\bar{s}|<\eta_0,~\forall t>t_0$, which is a contradiction to $|s(t_k)-\bar{s}|\ge \eta_0$. Thus we have $|s(t)-\bar{s}|\ge \eta_0,~t>T$. Let
    $$\eta =\min \left\{\min_{0\le t\le T}s(t),\eta_0\right\}>0,$$
    then it holds that $|s(t)-\bar{s}|\ge \eta, ~\forall t >0$. Let
    $$D_{\varepsilon}=\{s\in \mathbb{R}^n_{\ge 0}~|~dist(s,\partial \mathbb{R}^n_{\ge 0}) \ge \varepsilon \textrm{ and } |s-\bar{s}|\ge \eta\},
    $$
    where $\varepsilon>0$, $dist(s,\partial \mathbb{R}^n_{\ge 0})=\inf_{x \in \partial \mathbb{R}^n_{\ge 0}}|x-s|$ represents the distance from $s$ to the boundary of $\mathbb{R}^n_{\ge 0}$. For $s \in D_\varepsilon$, it satisfies $\alpha(|s-\bar{s}|) \ge \alpha(\eta)$. Since $\rho(|u(t)-\bar{u}|) \to 0$, there exist $\delta >0,~T_1>0$ such that
    $$dV/dt \le -\alpha(|s(t)-\bar{s}|)+\rho(|u(t)-\bar{u}|)\le -\delta,~t>T_1.
    $$
    Hence, when $t>T_1$, the amount of time that any trajectory spends in the set $D_\varepsilon$ is bounded above by $V(s(T_1))/ \delta$ (otherwise we have $s(t) \to \bar{s}$). That is to say, there exists a $T_2>T_1>0$ such that $dist\left(s(t),\partial \mathbb{R}^n_{\ge 0}\right)<\varepsilon,~t>T_2$. Since $\varepsilon$ is arbitrary, we have $s(t) \to \partial \mathbb{R}^n_{\ge 0}$. $\hfill\square$

\textbf{Proof of Lemma \ref{lem2}:} The proof is straightforward by combining \textit{Lemma} \ref{lem1} and the definition of persistence. $\hfill\square$

\textbf{Proof of Theorem \ref{thm_main}:} We just need to prove that the solution to the system (\ref{coupled_system}) satisfies $y^{\mathtt{2}}(t) \to \bar{y}^{\mathtt{2}}$ as $t \to +\infty$. By treating $y^{\mathtt{2}}$ in (\ref{coupled_system}) as the state variable and considering the remaining variables as the input $u(t)$, (\ref{coupled_system}) can be reformulated in the standard form 
    \begin{equation}
        \dot{y}^{\mathtt{2}}=f(u,y^{\mathtt{2}}), \quad y^{\mathtt{2}}(0)=y^{\mathtt{2}}_0,
    \end{equation}
    which is the dynamics of $\tilde{\mathscr{C}}^{\mathtt{2}}$. Since $\tilde{\mathscr{C}}^{\mathtt{2}}$ satisfies condition (1)-(3), according to \textit{Theorem 3.6} in \citep{chaves2005input}, there is an ISS-Lyapunov function $V(y^{\mathtt{2}})$ that satisfies condition (1)-(4) in \textit{Lemma} \ref{lem1}. Condition (3) also implies that $y^{\mathtt{2}}$ is bounded. In addition, note that $\tilde{\mathscr{C}}^{\mathtt{2}}$ has single linkage class, by \textit{Corollary 4.12} in \citep{anderson2011proof}, the trajectory $y^{\mathtt{2}}(t)$ is persistent. Therefore, by \textit{Lemma} \ref{lem2}, we have $y^{\mathtt{2}} \to \bar{y}^{\mathtt{2}}$.$\hfill\square$

\end{document}